\documentclass[reqno,12pt]{amsart}
\usepackage{latexsym,amssymb,amsfonts, amsmath}
\usepackage{a4wide}
\newtheorem{theorem}{Theorem}[section]
\newtheorem{corollary}[theorem]{Corollary}
\newtheorem{lemma}[theorem]{Lemma}
\newtheorem{proposition}[theorem]{Proposition}
\theoremstyle{definition}

\newtheorem{remarks}[theorem]{Remarks}
\newtheorem{example}[theorem]{Example}
\numberwithin{equation}{section}

\begin{document}
\bibliographystyle{acm}

\author{Gordon Blower}
\address{Department of Mathematics and Statistics, Lancaster University, Lancaster LA1 4YF, UK}
\email{g.blower@lancaster.ac.uk} 

\author{Fran\c cois Bolley}
\address{ENS Lyon, Umpa, 46 all\'ee d'Italie, F-69364 Lyon Cedex 07}
\email{fbolley@umpa.ens-lyon.fr}

\title[]
{Concentration of measure on product spaces with applications to Markov 
processes}

\thanks{2000 Mathematics Subject Classification : 60E15 (60E05 39B62)} 

\thanks{16th May 2005}

\keywords{Logarithmic Sobolev inequality, optimal transportation}

\begin{abstract}
For a stochastic process with state space some Polish space, 
this paper gives sufficient conditions on the initial and conditional distributions for the 
joint law to satisfy Gaussian concentration inequalities
and transportation inequalities. In the case of the Euclidean space ${\Bbb R}^{m}$, there are
sufficient conditions for the joint law to satisfy a logarithmic Sobolev inequality. 
In several cases, the obtained
constants are of optimal order of growth with respect to the number of random 
variables, or are independent of this number. These results extend results known for
mutually independent random variables to weakly dependent random variables under
Dobrushin--Shlosman type conditions. The paper also contains applications to 
Markov processes including the ARMA process. 
\end{abstract}

\maketitle

\section{Introduction} 

Given a complete and separable metric space $(X,d)$, ${\hbox{Prob}}(X)$ 
denotes the space of Radon probability measures on $X$, equipped with the (narrow) weak topology.
We say that $\mu\in {\hbox{Prob}}\, (X$) satisfies a 
{\sl Gaussian concentration inequality} $GC(\kappa)$ with constant $\kappa$ on $(X,d)$ if 
$$
\int_{X}\exp\bigl({tF(x)}\bigr) \mu (dx)\leq \exp\Bigl({t\int_{X} F(x)\, \mu(dx)
 +\kappa \, t^2/2}\Bigr)\qquad (t\in {\Bbb R})$$
holds for all $1$-Lipschitz functions $F:(X,d)\rightarrow {\Bbb R}$ (see \cite{BG99}). Recall that a function $g:(\Omega_1, d_1)\rightarrow (\Omega_2, d_2)$
between metric spaces is $L$-{\sl Lipschitz} if $d_2(g(x), g(y))\leq L \, d_1(x,y)$ holds for 
all $x,y\in \Omega_1$, and we call the infimum of such $L$ the
{\sl Lipschitz seminorm} of $g.$

\medskip

For $k \geq1$ and $x_1, \dots, x_k$ in $X$, we let $x^{(k)} = 
(x_1, \dots, x_k) \in X^k$ and, given 
$1\leq s<\infty$, we equip the product space $X^k$ with the metric $d^{(s)}$ defined
by $d^{(s)}(x^{(k)},y^{(k)})=(\sum_{j=1}^k d(x_j,y_j)^s)^{1/s}$ for $x^{(k)}$ and $y^{(k)}$ in $X^k$.\par

\indent Now let $(\xi_j )_{j=1}^n$ be a stochastic process with state space $X$. The first aim of this paper is 
to obtain concentration inequalities for the joint distribution  $P^{(n)}$ of 
$\xi^{(n)}=(\xi_1,\dots ,\xi_n)$, under hypotheses on the initial distribution $P^{(1)}$ of  $\xi_1$ and 
the conditional distributions $p_k(. \mid x^{(k-1)})$ of $\xi_k$  given $\xi^{(k-1)}$;
we recall that $P^{(n)}$ is given by
$$
P^{(n)} (dx^{(n)}) = p_n(dx_n \mid x^{(n-1)}) \dots \, p_2(dx_2 \mid x_1) \, P^{(1)}(dx_1).
$$
\noindent If the $(\xi_j)_{j=1}^n$ are mutually independent, and the
distribution of each $\xi_j$ satisfies $GC(\kappa )$, then 
$P^{(n)}$ on $(X^n, d^{(1)})$ is the product of the marginal 
distributions, and inherits $GC(n\kappa )$ from its marginal distributions by a 
simple `tensorization' argument. A similar result also applies to product measures for 
the transportation and logarithmic Sobolev inequalities which we
consider later; see \cite{Led01, Tal96}. To obtain concentration inequalities for $P^{(n)}$ when
 $(\xi_j)$ are weakly dependent, we impose additional
restrictions on the coupling between the variables, expressed in terms of Wasserstein distances which are 
defined as follows.

\medskip

Given $ 1 \leq s < \infty$, ${\hbox{Prob}}_s(X)$ denotes the
subspace of ${\hbox{Prob}}(X)$ consisting of $\nu$ such that $\int_X d(x_0,y)^s \, \nu (dy)$ is finite 
for some or equivalently all $x_0\in X$.
Then we define the {\sl Wasserstein distance} of order $s$ between $\mu$ and $\nu$ in ${\hbox{Prob}}_s(X)$ 
by
\begin{equation}\label{1.1}
W_s(\mu ,\nu )=\inf_{\pi} 
\Bigl( \int\!\!\!\int_{X\times X} d(x,y)^s \, \pi (dx \, dy)\Bigr)^{1/s}
\end{equation}
where $\pi\in {\hbox{Prob}}_s\, (X \times X)$ has marginals $\pi_1=\mu$ and 
$\pi_2=\nu$. Then $W_s$ defines a metric on ${\hbox{Prob}}_s(X)$, which in turn 
becomes a complete and separable metric space (see \cite{Rac91, Vil03}).

\smallskip

In section 3 we obtain the following result for time-homogeneous 
Markov chains.

\begin{theorem}\label{thm1.1}
Let $(\xi_j)_{j=1}^n$ be an homogeneous Markov process with
state space $X$, initial distribution $P^{(1)}$ and transition measure $p(. \mid x)$. Suppose that there exist
constants $\kappa_1$ and $L$ such that:

(i) $P^{(1)}$ and $p(. \mid x) $ ($x \in X$) satisfy $GC(\kappa_1)$ on $(X,d)$;

(ii) $x \mapsto p(. \mid x)$ is $L$-Lipschitz $(X,d) \to
({\hbox{Prob}}_1\,(X), W_1)$.

\noindent Then the joint law $P^{(n)}$ of $(\xi_1, \dots, \xi_n)$ satisfies $GC(\kappa_n)$ on
 $(X^n, d^{(1)})$, where
$$
\kappa_n = \kappa_1 \, \sum_{m=1}^{n} \Bigl( \sum_{k=0}^{m-1} L^k \Bigr)^2.
$$
\end{theorem}

In Example \ref{ex6.3} we demonstrate sharpness of these constants by providing 
for each value of $L$ a process such that $\kappa_n$ has optimal
growth in $n$.

\bigskip

Concentration inequalities are an instance of the wider class of
transportation inequalities, which bound the transportation cost by the relative entropy.
We recall the definitions. 

\smallskip

Let $\nu$ and $\mu$ be in ${\hbox{Prob}}(X)$, where $\nu$ is absolutely
continuous with respect to $\mu$, and let $d\nu /d\mu$ be the Radon--Nikodym
derivative. Then we define the {\sl relative entropy} of $\nu$ with respect to $\mu$ by
$$
{\hbox{Ent}}(\nu \mid \mu )=\int_X \log {{d\nu}\over {d\mu}}\, d\nu \, ;
$$
\noindent note that $0\leq {\hbox{Ent}}(\nu \mid \mu )\leq \infty$ by Jensen's inequality. By
convention we let ${\hbox{Ent}}(\nu \mid \mu ) =\infty$ if $\nu$ is not absolutely
continuous with respect to $\mu$.\par
\indent Given $1 \leq s < \infty$, we say that $\mu\in {\hbox{Prob}}_s(X)$ satisfies a 
{\sl transportation inequality} $T_s(\alpha )$ for cost function $d(x,y)^s$, with constant $\alpha$, if 
$$
W_s(\nu ,\mu )\leq \Bigl( {{2}\over {\alpha}}{\hbox{Ent}}(\nu \mid \mu)\Bigr)^{1/2}
$$
\noindent for all $\nu\in {\hbox{Prob}}_s(X)$.

Marton \cite{Mar96AP} introduced $T_2$ as `distance-divergence' 
inequalities in the context of information theory; subsequently Talagrand \cite{Tal96} showed that the 
standard Gaussian distribution on ${\Bbb R}^m$ satisfies $T_2(1).$ 
Bobkov and G\"otze  showed in \cite{BG99} that $GC(\kappa)$ is equivalent 
to $T_1(1/\kappa )$; their proof used the Kantorovich--Rubinstein duality
result, that
$$
W_1(\mu , \nu )=\sup_f \, \Bigl\{\int_X f(x)\mu (dx) -\int_X f(y) \, \nu(dy)\Bigr\} 
$$
where $\mu ,\nu\in{\hbox{Prob}}_1(X)$ and $f$ runs over the set of $1$-Lipschitz 
functions $f:X\rightarrow {\Bbb R}$. A $\nu\in {\hbox{Prob}}(X)$ satisfies a
$T_1$ inequality if and only if $\nu$ admits a square-exponential moment; that is,
$\int_X\exp({\beta d(x,y)^2})\,\nu (dx)$ is finite for some $\beta >0$ and some, and thus
all, $y\in X$; see \cite{BV05, DGW04}
 for detailed statements. Moreover, since $T_s(\alpha )$ implies $T_r(\alpha )$ for 
$1\leq r\leq s$ by H\"older's inequality, transportation inequalities are a tool for proving and 
strengthening concentration inequalities; they are also related to the Gaussian isoperimetric inequality 
as in \cite{Blo03}. For applications to empirical distributions in statistics, see \cite{Mas03}.

\medskip

Returning to weakly dependent $(\xi_j )_{j=1}^n$ with state space $X$,
we obtain transportation inequalities for the joint distribution $P^{(n)}$, under 
hypotheses on $P^{(1)}$ and the conditional distributions. 
 Djellout, Guillin and 
Wu  \cite{DGW04} developed Marton's coupling method \cite{Mar96AP, Mar04} to prove $T_s(\alpha )$
for $P^{(n)}$ under 
various mixing or contractivity conditions; see also \cite{Sam00}, or \cite{BV05} where the conditions are 
expressed solely in terms of exponential moments. We extend these 
results in sections \ref{sec2} and \ref{sec3} below, thus obtaining a strengthened dual form of Theorem \ref{thm1.1}.

\begin{theorem}\label{thm1.2}
Let $(\xi_j)_{j=1}^n$ be an homogeneous Markov process with
state space $X$, initial distribution $P^{(1)}$ and transition measure $p(. \mid x)$. Suppose that there exist constants 
$1\leq s\leq 2$, $\alpha >0$ and $L\geq 0$ such that:

(i) $P^{(1)}$ and $p(.\mid x)$ ($x \in X$) satisfy $T_s(\alpha )$;

(ii) $x\mapsto p(. \mid x)$ is $L$-Lipschitz 
$(X,d)\rightarrow ({\hbox{Prob}}_s\, (X), W_s)$.

\noindent Then the joint distribution $P^{(n)}$ of $(\xi_1, \dots, \xi_n)$ 
satisfies $T_s(\alpha_n)$, where
$$
\alpha_n = \left\{ \begin{array}{ll}
		n^{1-(2/s)}(1-L^{1/s})^2 \alpha     & \qquad \mbox{if} \quad L<1,\\
		e^{(2/s)-2} (n^{-(2/s)-1} \alpha & \qquad \mbox{if} \quad L=1,\\  
		\displaystyle \Bigl(\frac{L-1}{e^{s-1}L^n}\Bigr)^{2/s} \frac{\alpha}{n+1}
						    &  \qquad \mbox{if} \quad L>1;
			\end{array}
			\right.			    
$$
in particular $\alpha_n$ is independent of $n$ for $s=2$ when $L<1$.
\end{theorem}

Our general transportation Theorem \ref{thm2.1} will involve processes that are not necessarily Markovian,
but satisfy some a hypothesis related to Dobrushin--Shlosman's mixing condition 
[8, p. 352; 15, Definition 2]. When $X={\Bbb R}^m$, we shall also present some more 
computable version of hypothesis (ii) in Proposition \ref{prop2.2}, and later consider a 
stronger functional inequality.

\bigskip

A probability measure $\mu$ on ${\Bbb R}^m$ satisfies the {\sl logarithmic Sobolev 
inequality} $LSI(\alpha)$ with constant $\alpha>0$ if
$$
\int_{{\Bbb R}^m} f^2 \log \Bigl( f^2\slash \int_{{\Bbb R}^m}f^2 \, d\mu \Bigr) d\mu 
\leq {{2}\over{\alpha}}\int_{{\Bbb R}^m} \Vert \nabla f\Vert_{{\ell^2}}^2 \, d\mu 
$$
\noindent holds for all $f\in L^2(d\mu )$ that have distributional
gradient $\nabla f \in L^2(d\mu; {\Bbb R}^m)$. Given 
$(a_k)\in {\Bbb R}^m$, let $\Vert (a_k)\Vert_{\ell^s}=(\sum_{k=1}^{m}\vert a_k\vert^s)^{1/s}$ 
for $1 \leq s < \infty$, and 
$\Vert (a_k)\Vert_{\ell^{\infty}} = \sup_{1 \leq k \leq m} \vert a_k\vert$.

The connection
between the various inequalities is summarized by
\begin{equation}\label{1.2}
LSI(\alpha )\Rightarrow T_2(\alpha )\Rightarrow T_1(\alpha )\Leftrightarrow
GC(1/\alpha );
\end{equation}
see [3; 18; 24, p. 293]. Conversely, Otto and Villani showed that if
 $\mu (dx) = e^{-V(x)}\, dx$ satisfies 
$T_2(\alpha )$ where $V:{\Bbb R}^m\rightarrow {\Bbb R}$
is convex, then $\mu$ also satisfies $LSI(\alpha /4)$ (see [4; 18; 24, p. 298]); but this
converse is not generally true, as a counter-example in \cite{CG05} shows.

Gross  \cite{Gro75} proved that the standard Gaussian probability measure on 
${\Bbb R}^m$ satisfies $LSI(1)$. More generally, Bakry and Emery  \cite{BE85}
showed that if $V$ is twice continuously differentiable, with ${\hbox{Hess}}\, V\geq \alpha I_m$ on ${\Bbb
R}^m$ for some $\alpha >0$, then $\mu(dx) = e^{-V(x)}\, dx$ satisfies
$LSI(\alpha)$; see for instance \cite{Wan97} for extensions to this result. 
Whereas Bobkov and G\"otze \cite{BG99} characterized in
terms of their cumulative distribution functions those $\mu\in {\hbox{Prob}}
({\Bbb R})$ that satisfy $LSI(\alpha )$ for some $\alpha$, there is no known geometrical 
characterization of such probability measures on ${\Bbb R}^m$ when $m>1$.

\medskip

Our main Theorem \ref{thm5.1} gives a sufficient condition for the joint law of a weakly
dependent process with state space ${\Bbb R}^m$ to 
satisfy $LSI$. In section \ref{sec6} we deduce the following for distributions of 
time-homogeneous Markov processes. Let $\partial /\partial x$ denote the gradient
with respect to $x\in {\Bbb R}^m$.

\begin{theorem}\label{thm1.3}
Let $(\xi_j)_{j=1}^n$ be an homogeneous Markov process with
state space ${\Bbb R}^m$, initial distribution $P^{(1)}$ and transition measure $p(dy \mid x) =e^{-u(x,y)}dy$. 
Suppose that there exist constants $\alpha >0$ and $L\geq 0$ such that:

(i) $P^{(1)}$ and $p(.\mid x)$ $(x\in {\Bbb R}^m)$ satisfy $LSI(\alpha )$;

(ii) $u$ is twice continuously differentiable and the off-diagonal blocks of 
its Hessian matrix satisfy 
$$
\Bigl\Vert {{\partial^2u}\over {\partial x\partial y}}\Bigr\Vert \leq L
$$
as operators $({\Bbb R}^m, \ell^2)\rightarrow ({\Bbb R}^m, \ell^2)$.

\noindent Then the joint law $P^{(n)}$ of the first $n$ variables $(\xi_1, \dots, \xi_n)$ satisfies $LSI(\alpha_n)$,
 where
$$
\alpha_n = \left\{ \begin{array}{ll}
	\displaystyle	\frac{(\alpha - L )^2}{\alpha}     & \qquad \mbox{if} \quad L < \alpha,\\
	\displaystyle   \frac{\alpha}{n(n+1)(e-1)}         & \qquad \mbox{if} \quad L = \alpha,\\  
	\displaystyle	\Bigl( \frac{\alpha}{L} \Bigr)^{2n} \frac{L^2 - \alpha^2}{\alpha e(n+1)}
						   & \qquad \mbox{if} \quad  L > \alpha;
			\end{array}
			\right.			    
$$
in particular $\alpha_n$ is independent of $n$ when $L<\alpha$.
\end{theorem}

\medskip

The plan of the paper is as follows. In section \ref{sec2} we state and prove our 
results on transportation inequalities, which imply Theorem \ref{thm1.2}, and in section
\ref{sec3} we deduce Theorem \ref{thm1.1}. In section \ref{sec4} we prove $LSI(\alpha )$ for the joint
distribution of 
ARMA processes, with $\alpha$ independent of the size of the sample. 
In section \ref{sec5} we obtain a more
general $LSI$, which we express in a simplified form for Markov processes in section \ref{sec6}. 
Explicit examples in section \ref{sec6} show that several of our results have optimal growth of 
the constants with respect to $n$ as $n\rightarrow\infty$, and that the 
hypotheses are computable and 
realistic.

\medskip

\section{Transportation inequalities}\label{sec2} 
Let $(\xi_k)_{k=1}^n$ be a stochastic process with state space $X$, let 
$p_k(. \mid x^{(k-1)})$ denote the transition measure 
between the states at times $k-1$ and $k$, and let $P^{(n)}$ be the joint
distribution of $\xi^{(n)}$. Our main result of this section is a transportation
inequality.

\begin{theorem}\label{thm2.1}
Let $1 \leq s \leq 2$, and suppose that there exist $\alpha_1>0$ 
and $M \geq \rho_{\ell} \geq 0$ $(\ell =1, \dots, n)$ such that:

(i) $P^{(1)}$ and $p_k(. \mid x^{(k-1)})$ $(k=2, \dots ,n; \, x^{(k-1)} \in X^{k-1})$ satisfy 
$T_s(\alpha )$ on $(X,d)$;

(ii) $x^{(k-1)}\mapsto p_{k}(.\mid x^{(k-1)})$ is Lipschitz as a map
$(X^{k-1}, d^{(s)})\rightarrow ({\hbox{Prob}}_s\, (X), W_s)$
for $k=2, \dots, n$, in the sense that
$$
W_s\bigl(p_k(.\mid x^{(k-1)}), p_k(.\mid y^{(k-1)})\bigr)^s\leq
\sum_{j=1}^{k-1}\rho_{k-j} \, d(x_j,y_j)^s\qquad (x^{(k-1)}, y^{(k-1)}\in X^{k-1}).
$$

Then $P^{(n)}$ satisfies the transportation inequality $T_s(\alpha_n)$ where
$$
\alpha_n =\alpha \left({{(n \, e)^{1-s} M}\over {(1+M)^{n}}}\right)^{2/s}.
$$

Suppose further that 

(iii) $\sum_{j=1}^n\rho_j \leq R$.

\noindent Then the joint distribution $P^{(n)}$ satisfies $T_s(\alpha_n)$ where
$$
\alpha_n = \left\{ \begin{array}{ll}
		 n^{1-(2/s)}(1-R^{1/s})^2\alpha  & \qquad \mbox{if} \quad R < 1,\\
		 e^{(2/s)-2} (n+1)^{-(2/s)-1}    & \qquad \mbox{if} \quad R = 1,\\  
	\displaystyle	\Bigl(\frac{R-1}{e^{s-1}R^n}\Bigr)^{2/s} \frac{\alpha}{n+1}
						 & \qquad \mbox{if} \quad R > 1.
			\end{array}
			\right.			    
$$
\end{theorem}

\medskip

In hypothesis (iii), the sequence $(\rho_k)_{k=1}^{n-1}$ measures the extent to
which the distribution of $\xi_n$ depends upon the previous $\xi_{n-1}, \xi_{n-2}, \dots $; so in
most examples $(\rho_k)_{k=1}^{n-1}$ is decreasing. \par
A version of Theorem \ref{thm2.1} was obtained by Djellout, Guillin and Wu, but
with an explicit constant only when $R<1$; see \cite[Theorem 2.5 and Remark 2.9]{DGW04}. Theorem 
\ref{thm2.1} also improves upon section 4 of \cite{BV05}, where the 
assumptions were written in terms of moments of the considered measures.\par
\indent The Monge--Kantorovich transportation problem involves finding, for given $\mu,\nu \in
{\hbox{Prob}}(X)$, an optimal transportation strategy in 
\eqref{1.1}, namely a $\pi $ that minimises the transportation cost; a compactness and semi-continuity argument 
ensures that, for suitable cost functions, there always exists 
such a $\pi$. We recall that, given $\mu\in {\hbox{Prob}}(X)$, another Polish space $Y$ and a continuous 
function $\varphi :X\rightarrow Y$, the measure {\sl induced} from $\mu$ by $\varphi$ 
is the unique $\nu\in {\hbox{Prob}}(Y)$ such that 
$$
\int_Y f(y)\nu (dy)=\int_X f(\varphi (x))\mu (x)
$$
for all bounded and continuous $f:X\rightarrow {\Bbb R}$. 
Brenier and McCann showed that if $\mu$ and $\nu$ belong to ${\hbox{Prob}}_2({\Bbb R}^m)$, and if moreover 
$\mu$ is absolutely continuous with respect to Lebesgue measure, 
then there exists a convex function $\Phi :{\Bbb R}^m\rightarrow {\Bbb R}$ such that the gradient 
$\varphi=\nabla \Phi$ induces $\mu$ from $\nu$ and gives the unique solution to the Monge--Kantorovich
transportation problem for $s=2$, in the sense that
$$
\int_{{\Bbb R}^m} \Vert \nabla \Phi(x) - x \Vert^2_{\ell^2} \, \mu(dx) = W_2(\mu,\nu)^2.
$$
Further extensions of this result were obtained by Gangbo and McCann for
$1<s\leq 2$, by Ambrosio and Pratelli for $s=1$, and by McCann \cite{MC01} 
in the context of compact and connected $C^3$-smooth Riemannian manifolds that are without boundary 
(see also \cite{CMS01, Vil03}).

\smallskip

\begin{proof}[Proof of Theorem \ref{thm2.1}]
In order to give an explicit solution in a case
of importance, we first suppose that $X= {\Bbb R}^m$ and that $P^{(1)}$ and 
$p_{j}(dx_j \mid x^{(j-1)})$ ($j=2, \dots, n$) are all absolutely continuous with respect to Lebesgue 
measure. Then let $Q^{(n)}\in {\hbox{Prob}}_s ({\Bbb R}^{nm})$ be of finite relative entropy with respect to 
$P^{(n)}$. Let $Q^{(j)}(dx^{(j)})$ be the marginal distribution of $x^{(j)}\in {\Bbb R}^{jm}$ with respect to 
$Q^{(n)}(dx^{(n)})$, and disintegrate $Q^{(n)}$ in terms of 
conditional probabilities, according to
$$
Q^{(j)}(dx^{(j)})=q_j(dx_j\mid x^{(j-1)})Q^{(j-1)}(dx^{(j-1)}).
$$
In particular $q_j(. \mid x^{(j-1)})$ is absolutely continuous with
respect to $p_j(. \mid x^{(j-1)})$ and hence with respect to Lebesgue measure, for $Q^{(j-1)}$ 
almost every $x^{(j-1)}$. A standard computation ensures that

\begin{eqnarray}\label{2.1}
{\hbox{Ent}}\, (Q^{(n)}\mid P^{(n)})
& = &  {\hbox{Ent}}\, (Q^{(1)}\mid P^{(1)}) \\
&   & \displaystyle +\sum_{j=2}^n\int_{{\Bbb R}^{(j-1)m}} {\hbox{Ent}}\, \bigl(q_j(\,.\,\mid x^{(j-1)})\mid
p_j(\,.\,\mid x^{(j-1)})\bigr)\,Q^{(j-1)}(dx^{(j-1)}) . \nonumber
\end{eqnarray}

When the hypothesis (i) of Theorem 2.1 holds for 
some $1<s\leq 2$, it also holds for $s=1$. 
Consequently, by the Bobkov--G\"otze theorem, $P^{(1)}$ and 
$p_j(dx_j\mid x^{(j-1)})$ satisfy 
$GC(\kappa )$ for $\kappa =1/\alpha $, and then one can check that there exists $\varepsilon >0$
such that 
$$
\int_{{\Bbb R}^{m}}\exp (\varepsilon \Vert x^{(1)} \Vert_{\ell^2}^2)P^{(1)}(dx^{(1)})<\infty
$$
and likewise for $p_j$; compare with Herbst's theorem \cite[p. 280]{Vil03}, and \cite{BG99, DGW04}.
Hence $Q^{(1)}$ and $q_j(dx_j\mid x^{(j-1)})$ for $Q^{(j-1)}$ almost every $x^{(j-1)}$ have finite 
second moments, since by Young's inequality
$$
\int_{{\Bbb R}^{m}} \varepsilon \Vert x^{(n)}\Vert^2_{\ell^2} \, Q^{(1)}(dx^{(1)})
\leq {\hbox{Ent}}(Q^{(1)}\mid P^{(1)})+
\log \int_{{\Bbb R}^{m}}\exp\bigl( \varepsilon \Vert x^{(1)}\Vert^2_{\ell^2}
\bigr)P^{(1)}(dx^{(1)})<\infty 
$$
\noindent and likewise with $q_j$ and $p_j$ in place of $Q^{(1)}$ and $P^{(1)}$
respectively.

Let $\theta_1:{\Bbb R}^m\rightarrow {\Bbb R}^m$ be an optimal
transportation map that induces $P^{(1)}(dx_1)$ from
$Q^{(1)}(dx_1)$; then for $Q^{(1)}$ every each $x_1$, let $x_2\mapsto \theta_2(x_1,x_2)$ induce
$p_2(dx_2\mid \theta_1(x_1))$ from $q_2(dx_2\mid x_1)$ optimally; hence 
$\Theta^{(2)}:{\Bbb R}^{2m}\rightarrow {\Bbb R}^{2m}$, defined by\par
\noindent $\Theta^{(2)}(x_1,x_2)=(\theta_1(x_1),\theta_2(x_1,x_2))$ on a certain set 
of full $Q^{(2)}$ measure, induces $P^{(2)}$ from $Q^{(2)}.$
Generally, having constructed $\Theta^{(j)}:{\Bbb R}^{jm}\rightarrow {\Bbb R}^{jm}$, we
let $x_{j+1}\mapsto \theta_{j+1}(x^{(j)},x_{j+1})$ be an optimal transportation map that induces
$p_{j+1}(dx_{j+1}\mid \Theta^{(j)}(x^{(j)}))$ from $q_{j+1}(dx_{j+1}\mid x^{(j)})$, for 
all $x^{(j)}$ in 
a certain set of full $Q^{(j)}$ measure; then
we let $\Theta^{(j+1)}:{\Bbb R}^{(j+1)m}\rightarrow  {\Bbb R}^{jm}\times {\Bbb R}^{m}$
be the map defined by
$$
\Theta^{(j+1)}(x^{(j+1)})=(\Theta^{(j)}(x^{(j)}), \theta_{j+1}(x^{(j+1)}))$$
\noindent on a set of full $Q^{(j+1)}$ measure. In particular
$\Theta^{(j+1)}$ induces $P^{(j+1)}$ from $Q^{(j+1)}$, in the style of Kneser.

\smallskip

This transportation strategy may not be optimal, nevertheless it gives the bound
\begin{equation}\label{2.2}
W_s(Q^{(n)}, P^{(n)})^s\leq \int_{{\Bbb R}^{nm}}\Vert
\Theta^{(n)}(x^{(n)})-x^{(n)}\Vert^s_{\ell^s} Q^{(n)}(dx^{(n)}) = \sum_{k=1}^{n} d_k
\end{equation}
by the recursive definition of $\Theta^{(n)},$ where we have let
$$
d_k=\int_{{\Bbb R}^{km}}\Vert \theta_k(x^{(k)})-x_k\Vert^s_{\ell^s}
Q^{(k)}(dx^{(k)}) \qquad (k = 1, \dots, n).$$
\noindent However, the transportation at step $k$ is optimal by construction, so 
\begin{equation}\label{2.3}
d_k = \int_{{\Bbb R}^{(k-1)m}}
W_s\bigl( p_k(\,.\, \mid \Theta^{(k-1)}(x^{(k-1)})), q_k(\,.\,\mid
x^{(k-1)})\bigr)^sQ^{(k-1)}(dx^{(k-1)}). 
\end{equation}

Given $a,b>0$, $1\leq s\leq 2$ and $\gamma >1$, we have $(a+b)^s\leq
(\gamma/(\gamma -1))^{s-1}a^s+\gamma^{s-1}b^s$. Hence by the triangle inequality, the
expression \eqref{2.3} is bounded by
\begin{multline}\label{2.4}
\Bigl({{\gamma}\over {\gamma -1}}\Bigr)^{s-1}
\int_{{\Bbb R}^{m(k-1)}} W_s\bigl(p_k(\,.\,\mid x^{(k-1)}), q_k(\,.\,\mid x^{(k-1)})\bigr)^s
Q^{(k-1)}(dx^{(k-1)}) \\ 
\quad +\gamma^{s-1}\int_{{\Bbb R}^{m(k-1)}}W_s( p_k
(\,.\, \mid \Theta^{(k-1)}(x^{(k-1)})), p_k(\,.\,\mid x^{(k-1)}))^s 
Q^{(k-1)}(dx^{(k-1)}).
\end{multline}
By hypothesis (i) and then H\"older's inequality, we bound the first integral in \eqref{2.4} by 
$$
h_k = \Bigl({{\gamma}\over {\gamma -1}}\Bigr)^{s-1} \Bigl({{2}\over {\alpha}}\Bigr)^{s/2}  
\Bigl( \int_{{\Bbb R}^{(k-1)m}}{\hbox{Ent}}(q_k\mid p_k) \, dQ^{(k-1)}\Bigr)^{s/2}.
$$
Meanwhile, on account of hypothesis (ii) the second integral in \eqref{2.4} is bounded by
$$
\gamma^{s-1} \int_{{\Bbb R}^{m(k-1)}} \sum_{j=1}^{k-1} \rho_{k-j} \bigl\Vert
\theta_j(x^{(j)})-x_j\bigr\Vert^s \, Q^{(k-1)}(dx^{(k-1)})
= \gamma^{s-1} \sum_{j=1}^{k-1} \rho_{k-j} d_j,
$$
and when we combine these contributions to \eqref{2.4} we have 
\begin{equation}\label{2.5}
d_k \leq h_k + \gamma^{s-1} \sum_{j=1}^{k-1} \rho_{k-j} \, d_j. 
\end{equation}

In the case when the $\rho_{\ell}$ are merely bounded by $M$, one can prove by induction that
$$
d_k \leq h_k +  \gamma^{s-1} M \sum_{j=1}^{k-1} h_j (1 +  \gamma^{s-1} M)^{k-1-j},
$$
so that
$$
\sum_{k=1}^{n} d_k \leq \sum_{j=1}^{n} h_j (1 + \gamma^{s-1} M)^{n-j} \leq 
\Bigl(\sum_{j=1}^{n} h_j^{2/s} \Bigr)^{s/2} 
\Bigl(\sum_{j=1}^{n} (1+ \gamma^{s-1} M)^{2 (n-j)/{(2-s)}} \Bigr)^{({2-s})/{2}}
$$
by H\"older's inequality. The first sum on the right-hand side is
$$
\Bigl(\sum_{j=1}^{n} h_j^{2/s} \Bigr)^{s/2} = 
\Bigl({{\gamma}\over {\gamma -1}}\Bigr)^{s-1} \Bigl({{2}\over {\alpha}}\Bigr)^{s/2}
{\hbox{Ent}}\, (Q^{(n)}\mid P^{(n)})^{s/2}
$$
by \eqref{2.1}. Finally, setting $\gamma =1+1/n$, we obtain by \eqref{2.2} the stated result
$$
W_s(Q^{(n)}, P^{(n)})^s \leq
\Bigl({{2}\over {\alpha}}\Bigr)^{s/2} \frac{(1+M)^n}{M} \, (n \, e)^{s-1}
\, {\hbox{Ent}}(Q^{(n)}\mid P^{(n)})^{s/2}.
$$

\medskip

(iii) Invoking the further hypothesis (iii), we see that $T_m = \sum_{j=1}^{m} d_j$
satisfies on account of \eqref{2.5} the recurrence relation 
$$
T_{m+1} \leq \sum_{j=1}^{m+1} h_j + \gamma^{s-1} R \, T_m,
$$ 
which enables us to use H\"older's inequality again and bound $T_n$ by 
\begin{multline*}
\sum_{k=1}^{n} \bigl( \sum_{j=1}^{k} h_j \bigr)  (\gamma^{s-1} R)^{n-k}
= \sum_{j=1}^{n} h_j \sum_{\ell=0}^{n-j} (\gamma^{s-1} R)^{\ell} \\
\leq  \Bigl( \sum_{j=1}^{n} h_j^{2/s} \Bigr)^{s/2} 
     \Bigl( \sum_{j=1}^{n} \bigl(\sum_{\ell=0}^{n-j} 
(\gamma^{s-1} R)^{\ell} \bigr)^{{2}/({2-s})}
     \Bigr)^{{(2-s)}/{2}}
\end{multline*}
for $1\leq s< 2$. By \eqref{2.2} and the definition of $T_n$ this leads to
\begin{eqnarray}
W_s(Q^{(n)},P^{(n)})^s 
&\leq & \Bigl({{\gamma}\over {\gamma -1}}\Bigr)^{s-1}
\Bigl( \sum_{m=1}^{n} \bigl(\sum_{\ell=0}^{m-1} (\gamma^{s-1} R)^{\ell} \bigr)^{{2}/({2-s)}}
     \Bigr)^{{(2-s)}/{2}} 
\Bigl( {{2}\over{\alpha}} \, {\hbox{Ent}}(Q^{(n)}\mid P^{(n)})\Bigr)^{s/2} \label{2.6}\\
&\leq & \Bigl({{\gamma}\over {\gamma -1}}\Bigr)^{s-1} n^{1-s/2} \sum_{\ell=0}^{n-1} (\gamma^{s-1} R)^{\ell}
\, \Bigl( {{2}\over{\alpha}} \, {\hbox{Ent}}(Q^{(n)}\mid P^{(n)})\Bigr)^{s/2}; \label{2.7}
\end{eqnarray}
this also holds for $s=2.$ Finally we select $\gamma$ according to the value of
 $R$ to make the bound \eqref{2.7} precise. When $R<1$, we let $\gamma =R^{-1/s}>1$, so that 
$\gamma^{s-1}R=R^{1/s}<1$, and we deduce the transportation inequality
$$
W_s(Q^{(n)},P^{(n)})^s\leq\Bigl( {{2}\over{\alpha}}\Bigr)^{s/2}{{n^{1-s/2}}\over
{(1-R^{1/s})^{s}}}{\hbox{Ent}}(Q^{(n)}\mid
P^{(n)})^{s/2}.
$$
When $R\geq 1$, we let $\gamma =1+1/n$ to obtain the transportation inequality 
$$
W_s(Q^{(n)},P^{(n)})^s \leq \Bigl({{2}\over{\alpha}}\Bigr)^{s/2}
(n+1)^{s-1} \, n^{1-s/2} \, 
\biggl({{(1+1/n)^{n(s-1)} R^n -1}\over{(1+1/n)^{s-1}R-1}}\biggr) \, 
{\hbox{Ent}}(Q^{(n)}\mid P^{(n)})^{s/2},
$$
which leads to the stated result by simple analysis, and
completes the proof when $X={\Bbb R}^m.$

\bigskip

For typical Polish spaces $(X,d)$, we cannot rely on the existence of
optimal maps, but we can use a less explicit inductive approach to construct the 
transportation strategy, as in \cite{DGW04}. Given $j=1, \dots, n-1$, assume that 
$\pi^{(j)}\in {\hbox{Prob}}(X^{2j})$ has marginals $Q^{(j)}(dx^{(j)})$ and $P^{(j)}(dy^{(j)})$ 
and satisfies
$$
W_s(Q^{(j)}, P^{(j)})^s\leq\int_{X^{2j}} \sum_{k=1}^j
d(x_k,y_k)^s\pi^{(j)}(dx^{(j)}dy^{(j)}).
$$
Then, for each $(x^{(j)},y^{(j)})\in X^{2j},$ let $\sigma_{j+1}(\, . \, \mid x^{(j)},
y^{(j)})\in {\hbox{Prob}}(X^2)$ be an optimal transportation strategy that has marginals
$q_{j+1}(dx_{j+1}\mid x^{(j)})$ and $p_{j+1}(dy_{j+1}\mid y^{(j)})$ and that satisfies 
$$
 W_s(q_{j+1}(\,.\,\mid x^{(j)}), p_{j+1}(\,.\,\mid y^{(j)}))^s
=\int_{X^2}d(x_{j+1}, y_{j+1})^s \, \sigma_{j+1}(dx_{j+1} \, dy_{j+1}\mid x^{(j)}, y^{(j)}).$$
Now we let 
$$
\pi^{(j+1)}(dx^{(j+1)}dy^{(j+1)})=\sigma_{j+1}(dx_{j+1} \, dy_{j+1}\mid x^{(j)},y^{(j)})
\, \pi^{(j)}(dx^{(j)}dy^{(j)}),$$which defines a probability on $X^{2(j+1)}$ with marginals $Q^{(j+1)}(dx^{(j+1)})$
and $P^{(j+1)}(dy^{(j+1)})$. This may not give an optimal transportation strategy; 
nevertheless, the recursive definition shows that 
$$
W_s(Q^{(n)}, P^{(n)})^s \leq \sum_{j=1}^n \int_{X^{2(j-1)}}
W_s\bigl(q_j(\, .\,\mid x^{(j-1)}),p_j(\, .\,\mid
y^{(j-1)})\bigr)^s\pi^{(j-1)}(dx^{(j-1)}dy^{(j-1)})
$$
\noindent and one can follow the preceding proof from \eqref{2.2} onwards.
\end{proof}

\medskip

\begin{proof}[Proof of Theorem 1.2] 
Under the hypotheses of Theorem \ref{thm1.2}, we can take
$\rho_1=L$ and $\rho_j=0$ for $j=2, \dots, n$, which satisfy Theorem \ref{thm2.1} with
$R=L$ in assumption (iii).
\end{proof}

\medskip

The definition of $W_s$ not being well suited to direct calculation, 
we now give a computable sufficient condition for hypothesis (ii) of Theorem \ref{thm2.1} 
to hold with some constant coefficients $\rho_{\ell}$ when $(X,d)=({\Bbb R}^m, \ell^s)$.

\begin{proposition}\label{prop2.2}
Let $u_j:{\Bbb R}^{jm}\rightarrow {\Bbb R}$ be a twice
continuously differentiable function that has bounded second-order partial
derivatives. Let $1\leq s\leq 2$ and suppose further that:

(i) $p_j(dx_j\mid x^{(j-1)})=\exp( {-u_j(x^{(j)})})\, dx_j$ satisfies $T_s(\alpha)$ 
for some $\alpha >0$ and all $x^{(j-1)}\in {\Bbb R}^{m(j-1)}$;

(ii) there exists some real number $M_s$ such that 
$$
\sup_{x^{(j-1)}}\int_{{\Bbb R}^m}\Bigl\Vert
\Bigl(\frac{\partial u_j}{\partial x_k}\Bigr)_{k=1}^{j-1}\Bigr\Vert^2_{\ell^{s'}}p_j(dx_j\mid
x^{(j-1)})=M_s,
$$
where $1/s'+1/s=1$ and $\partial /\partial x_k$ denotes the gradient with
respect to $x_k$.

\noindent Then $x^{(j-1)}\mapsto p_j(. \mid x^{(j-1)})$ is $\sqrt{(M_s/
\alpha )}$-Lipschitz $({\Bbb R}^{m(j-1)},\ell^s)\rightarrow
({\hbox{Prob}}_s({\Bbb R}^m), W_s).$
\end{proposition}

\begin{proof} 
Given $x^{(j-1)}, \bar x^{(j-1)}\in {\Bbb
R}^{m(j-1)},$ we let $x^{(j-1)}(t)=(1-t)\bar x^{(j-1)}+tx^{(j-1)}$
$(0\leq t\leq 1)$ be the straight-line segment that joins them, and
we consider 
$$
f(t)=W_s\bigl(p_j(. \mid x^{(j-1)}(t)) , p_j(. \mid \bar x^{(j-1)})\bigr);
$$
then it suffices to show that $f:[0,1]\rightarrow {\Bbb R}$ is Lipschitz and to bound its 
Lipschitz seminorm.

By the triangle inequality and (i), we have
$$
\Bigl({{f(t+\delta )-f(t)}\over {\delta }}\Bigr)^2 
$$
$$ 
\leq  \frac{1}{\delta^2} W_s\bigl(p_j(. \mid x^{(j-1)}(t+\delta )), p_j(. \mid x^{(j-1)}(t))\bigr)^2 
$$
$$
\leq 
{{1}\over{\alpha\delta^2}}\Bigl\{ {\hbox{Ent}}\bigl( p_j(. \mid  x^{(j-1)}(t+\delta ))\mid p_j(. \mid x^{(j-1)}(t))\bigr)
 + {\hbox{Ent}}\bigl( p_j(. \mid 
x^{(j-1)}(t))\mid p_j(. \mid x^{(j-1)}(t+\delta ))\bigr)\Bigr\}
$$
\begin{multline}\label{2.8}
={{1}\over {\alpha\delta^2}}\int_{{\Bbb R}^m}\bigl(u_j(x^{(j-1)}(t+\delta ), x_j)-u_j(x^{(j-1)}(t),x_j)\bigr) \\
\bigl\{\exp (-u_j(x^{(j-1)}(t),x_j))-\exp (-u_j(x^{(j-1)}(t+\delta ),x_j))\bigr\}\,dx_j.
\end{multline}
However, by the assumptions on $u_j$ and the mean-value theorem, we have 
\begin{multline*}
 u_j(x^{(j-1)}(t+\delta),x_j)-u_j(x^{(j-1)}(t),x_j) \\
 =\delta \sum_{k=1}^{j-1}\bigl\langle \frac{\partial u_j}{\partial x_k} (x^{(j-1)}(t),x_j),
x_k-\bar x_k\bigr\rangle +{{\delta^2}\over {2}}\bigl\langle
{\hbox{Hess}}\,u_j\, (x^{(j-1)}-\bar x^{(j-1)}),(x^{(j-1)}-\bar
x^{(j-1)})\bigr\rangle,
\end{multline*}
where ${\hbox{Hess}}\, u_j$ is computed at some point
between $(x^{(j-1)},x_j)$ and $(\bar x^{(j-1)},x_j)$ and 
is uniformly bounded.  Proceeding in the same way for the other term \eqref{2.8}, we obtain
$$
\limsup_{\delta\rightarrow 0+}\Bigl( {{f(t+\delta
)-f(t)}\over {\delta }}\Bigr)^2\leq {{1}\over {\alpha}}\int_{{\Bbb
R}^m}\Bigl( \sum_{k=1}^{j-1}\bigl\langle \frac{\partial u_j}{\partial x_k} (x^{(j-1)}(t), x_j) , x_k-\bar
x_k\bigr\rangle\Bigr)^2 p_j(dx_j\mid x^{(j-1)}(t)).
$$
Hence by H\"older's inequality we have 
\begin{multline*}
\lim\sup_{\delta\rightarrow 0+}{{\vert f(t+\delta )-f(t)\vert}
\over {\delta}} \\
\leq {{1}\over {\sqrt{\alpha}}}\biggl(\int_{{\Bbb R}^m}
\Bigl(\sum_{k=1}^{j-1}\Bigl\vert \frac{\partial u_j}{\partial x_k} (x^{(j-1)}(t),x_j)\Bigr\vert^{s'}\Bigr)^{2/s'}
p_j(dx_j\mid x^{(j-1)}(t))\biggr)^{1/2}\Vert x^{(j-1)}-\bar x^{(j-1)}\Vert_{\ell^s}
\end{multline*}
for $1<s\leq 2$, and likewise with obvious changes for $s=1$.
By assumption (ii) and Vitali's theorem, $f$ is Lipschitz with constant 
$\sqrt{(M_s/ \alpha)}\Vert x^{(j-1)}-\bar
x^{(j-1)}\Vert_{\ell^s}$, as
required. 
\end{proof}

\medskip

\section{Concentration inequalities for weakly dependent sequences}\label{sec3} 
In terms of concentration inequalities, the dual version of Theorem \ref{thm2.1} 
reads as follows.

\begin{theorem}\label{thm3.1} 
Suppose that there exist $\kappa_1>0$ and $M\geq\rho_j\geq 0$ ($j=1, \dots, n$) such
that:

(i) $P^{(1)}$ and $p_k(. \mid x^{(k-1)})$ $(k=2, \dots ,n; \, x^{(k-1)} \in X^{k-1})$ satisfy 
$GC(\kappa_1 )$ on $(X,d)$;
 
(ii) $x^{(k-1)}\mapsto p_k(.\mid x^{(k-1)})$ is Lipschitz as a map 
$(X^{k-1}, d^{(1)}) \rightarrow ({\hbox{Prob}}_1\,(X), W_1)$ for $k=2, \dots, n$, in the sense that
$$
W_1\bigl(p_k(. \mid x^{(k-1)}), p_k(. \mid y^{(k-1)}) \bigr)\leq
\sum_{j=1}^{k-1} \rho_{k-j} \, d(x_j, y_j)\qquad (x^{(k-1)}, y^{(k-1)}\in X^{k-1}).
$$

\noindent Then the joint law $P^{(n)}$ satisfies $GC(\kappa_n)$ on $(X^n, d^{(1)})$, where
$$
\kappa_n =\kappa_1 {{(1+M)^{2n}} \over {M^2}}.
$$

Suppose moreover that 

(iii) $\sum_{j=1}^n\rho_j\leq R.$

\noindent Then $P^{(n)}$ satisfies $GC(\kappa_n (R))$ on $(X^n, d^{(1)})$, where 
$$
\kappa_n(R) = \kappa_1 \sum_{m=1}^{n} \Bigl( \sum_{k=0}^{m-1} R^k \Bigl)^2.
$$
\end{theorem}

\smallskip

\begin{proof}[Proof of Theorem 3.1] This follows from 
the Bobkov--G\"otze theorem \cite{BG99} and the bound \eqref{2.6} with $s=1$ 
in the proof of Theorem \ref{thm2.1}.
\end{proof}

\smallskip

Alternatively, one can prove Theorem \ref{thm3.1} directly by 
induction on the dimension, using the definition of $GC$.

\smallskip

\begin{proof}[Proof of Theorem 1.1] Under the hypotheses of Theorem \ref{thm1.1}, we can apply 
Theorem \ref{thm3.1} with $\rho_1 = L$ and $\rho_j = 0$ for $j=2, \dots, n$, which satisfy (iii)
\end{proof}

\medskip

\section{Logarithmic Sobolev inequalities for ARMA models}\label{sec4}
 In this section we 
give logarithmic Sobolev inequalities for the joint law
of the first $n$ variables from two auto-regressive moving average processes. 
In both results we obtain constants
that are independent of $n$, though the variables are not mutually independent, 
and we rely on the following general result which 
induces logarithmic Sobolev inequalities from one probability measure to another. 
For $m \geq1$,  let $\nu\in {\hbox{Prob}}({\Bbb R}^m)$ satisfy $LSI(\alpha)$, and
let $\varphi$ be a $L$-Lipschitz map from 
$({\Bbb R}^m, {\Bbb \ell}^{2})$ into itself; then, by the chain rule, the
probability measure that is induced from $\nu$ by $\varphi$ satisfies 
$LSI(\alpha / L^2)$. Our first application is the following.

\begin{proposition}\label{prop4.1} 
Let $Z_0$ and $Y_j$ $(j=1, 2, \dots )$ be 
mutually independent random variables in ${\Bbb R}^m$, and let 
$\alpha >0$ be a constant such that the distribution $P^{(0)}$ of $Z_0$ and the distribution
of $Y_j$ ($j =1, 2, \dots$) satisfy $LSI(\alpha)$.

\noindent Then for any $L$-Lipschitz map $\Theta :{\Bbb R}^m\rightarrow {\Bbb R}^m$,
the relation 
\begin{equation}\label{4.1}
Z_{j+1}=\Theta (Z_j)+Y_{j+1} \qquad (j=0,1, \dots )
\end{equation}
determines a stochastic process such that, for any $n \geq 1$, the joint distribution $P^{(n-1)}$ of
$(Z_j)_{j=0}^{n-1}$ satisfies $LSI(\alpha_n)$ where
$$
\alpha_n = \left\{ \begin{array}{ll}
			(1-L)^2 \alpha     & \qquad \mbox{if} \quad 0 \leq L < 1,\\
	\displaystyle   \frac{\alpha}{n(n+1)(e-1)}         & \qquad \mbox{if} \quad L = 1,\\  
	\displaystyle	\frac{L-1}{L^n} \frac{\alpha}{e(n+1)}
						   & \qquad \mbox{if} \quad  L > 1.
			\end{array}
			\right.			    
$$
\end{proposition}

\smallskip

\begin{proof}
For $(z_0, y_1, \dots, y_{n-1} )\in {\Bbb R}^{n m}$, let 
$\varphi_{n} (z_0, y_1, \dots, y_{n-1} )$ be the vector $(z_0, \dots, z_{n-1} )$,
defined by the recurrence relation 
\begin{equation}\label{4.2}
z_{k+1}=\Theta (z_k)+y_{k+1}\qquad (k=0, \dots, n-2 ).
\end{equation}
Using primes to indicate another solution of \eqref{4.2}, we deduce the following
inequality from the Lipschitz condition on $\Theta$: 
\begin{equation}\label{4.3}
\Vert z_{k+1}-z_{k+1}'\Vert^2\leq (1+\varepsilon)L^2 \Vert
z_k-z_k'\Vert^2+
(1+\varepsilon^{-1})\Vert y_{k+1}-y'_{k+1}\Vert^2
\end{equation} 
for all $\varepsilon >0$. In particular \eqref{4.3} implies the bound
$$
\Vert z_{k}-z_{k}'\Vert^2\leq \bigl((1+\varepsilon)L^2\bigr)^k \Vert z_0-z_0' \Vert^2 +
(1+\varepsilon^{-1}) \sum_{j=1}^{k} \bigl((1+\varepsilon)L^2\bigr)^{k-j} 
\Vert y_{j}-y'_{j}\Vert^2.
$$ 
By summing over $k$, one notes that $\varphi_{n}$ defines a Lipschitz function from 
$({\Bbb R}^{n m}, \ell^2)$ into itself, with Lipschitz seminorm 
$$
L_{\varphi_{n}}\leq\Bigl( {(1+\varepsilon^{-1})} \sum_{k=0}^{n-1} 
\bigl((1+\varepsilon) L^2\bigr)^k \Bigr)^{1/2}
$$
We now select $\varepsilon >0$ according to the value of $L$: when $L <1$, we let $\varepsilon = L^{-1} -1 >0$, 
so that $L_{\varphi_n} \leq (1-L)^{-1}$; whereas when $L \geq 1$, we let $\varepsilon = n^{-1}$, and obtain
$L_{\varphi_n} \leq [n(n+1)(e-1)]^{(1/2)}$ for $L=1$, and $L_{\varphi_n} \leq [e (n+1) L^n (L-1)^{-1}]^{1/2}$
for $L>1$.

Moreover, $\varphi_n$ induces the joint distribution of
$(Z_j)_{j=0}^{n-1} $ from the joint distribution of $(Z_0, Y_1, \dots, Y_{n-1})$. By independence, the joint 
distribution of $(Z_0, Y_1, \dots, Y_{n-1} )$ is a product measure on $({\Bbb R}^{n m}, \ell^2)$ that 
satisfies $LSI(\alpha )$. Hence the joint distribution of
$(Z_j)_{j=0}^{n-1}$ satisfies $LSI(\alpha )$, where $\alpha =L_{\varphi_n}^{-2}\, \alpha$.
\end{proof}

\medskip

\indent The linear case gives the following result for ARMA processes.

\begin{proposition}\label{prop4.2} 
Let $A$ and $B$ be $m\times m$ matrices such that the spectral radius
$\rho$ of $A$ satisfies $\rho <1$. Let also $Z_0$ and $Y_j$ $(j=1, 2, \dots )$ be mutually independent 
standard Gaussian $N(0,I_m)$ random variables in ${\Bbb R}^m$. Then, for any $n \geq 1$,  the joint distribution of 
the ARMA process $(Z_j)_{j=0}^{n-1}$, defined by 
the recurrence relation
$$
Z_{j+1}=AZ_j+BY_{j+1}\qquad (j=0, 1,\dots ),
$$
\noindent satisfies $LSI(\alpha )$ where 
$$
\alpha =\Bigl({{(1-\sqrt\rho )}\over {\max \{1,\Vert B\Vert\}}}\Bigr)^2
\Bigl( \sum_{j=0}^\infty \rho^{-j}\Vert A^j\Vert^2\Bigr)^{-2}.
$$
\end{proposition}

\smallskip

\begin{proof}

By Rota's Theorem \cite{Pau86}, $A$ is similar to a strict contraction on $({\Bbb
R}^m, \ell^2)$; that is, there
exists an invertible $m\times m$ matrix $S$ and a matrix $C$ such that $\Vert C\Vert\leq 1$ and 
$A=\sqrt \rho S^{-1}CS$; one can choose the similarity so that the
operator norms satisfy
$$
\Vert S\Vert \Vert S^{-1}\Vert \leq \sum_{j=0}^\infty \rho^{-j}\Vert A^j\Vert^2<\infty .
$$
Hence the ARMA process 
reduces to the solution of the recurrence relation
\begin{equation}\label{4.4}
SZ_{j+1}=\sqrt \rho \, CSZ_j+SBY_{j+1}\qquad (j=0, 1, \dots)
\end{equation}
which involves the $\sqrt\rho$-Lipschitz linear map $\Theta
:{\Bbb R}^m\rightarrow{\Bbb R}^m:$ $\Theta (w)=\sqrt \rho \, C \, w.$ Given
$n \geq 1$, the linear map $\Phi_n:{\Bbb R}^{n m}\rightarrow {\Bbb
R}^{n m}$, defined to solve \eqref{4.4} by
$$
(z_0, y_1, \dots ,y_n)\mapsto (Sz_0, SBy_1, \dots ,SBy_{n-1})\mapsto (Sz_0, Sz_1,\dots
,Sz_{n-1})\mapsto (z_0, z_1, \dots ,z_{n-1}),
$$
has operator norm 
$$
\Vert \Phi_n\Vert \leq \Vert S\Vert \Vert S^{-1}\Vert (1-\sqrt \rho )^{-1}\max \{
1, \Vert B\Vert \};
$$
moreover, $\Phi_n$ induces the joint distribution of $(Z_0, \dots ,Z_{n-1})$ from the joint
distribution of $(Z_0, Y_1, \dots ,Y_{n-1})$. By Gross's Theorem (see \cite{Gro75}), the latter distribution
satisfies $LSI(1)$, and hence the induced
distribution satisfies $LSI(\alpha )$, with $\alpha =\Vert
\Phi_n\Vert^{-2}$.
\end{proof}

\smallskip

\begin{remarks}\label{rem4.3} (i) As compared to Proposition \ref{prop4.1}, 
the condition imposed in Proposition \ref{prop4.2}
involves the spectral radius of the matrix $A$ and not its operator norm. In particular, for matrices
with norm $1$, Proposition \ref{prop4.1} only leads to $LSI$ with constant of order $n^{-2}$; whereas 
Proposition \ref{prop4.2} ensures $LSI$
with constant independent of $n$ under the spectral radius assumption $\rho <1$.

(ii) The joint distribution of the ARMA process is
discussed by Djellout, Guillin and Wu \cite[Section 3]{DGW04}. We have improved upon \cite{DGW04} by obtaining 
$LSI(\alpha)$, hence $T_2(\alpha )$, under the
spectral radius condition $\rho <1$, where $\alpha$ is independent of the size $n$ 
of the considered sample and the size of the matrices.
\end{remarks}

\medskip

\section{Logarithmic Sobolev inequality for weakly dependent processes}\label{sec5}
In this section we consider a stochastic process $(\xi_j)_{j=1}^n$, with state space ${\Bbb R}^m$ and
initial distribution $P^{(1)}$, which is not necessarily Markovian; we also assume that the transition 
kernels have positive densities with respect to Lebesgue measure, and write
$$
dp_j=p_j(dx_j\mid x^{(j-1)})=e^{-u_j(x^{(j)})} dx_j\qquad (j=2,\dots,n).
$$
The coupling between variables is measured by the following integral
$$
\Lambda_{j,k}(s) = \sup_{x^{(j-1)}}\int_{\Bbb R}\exp \Bigl( \bigl\langle s , {{\partial u_j}\over {\partial
x_k}}(x^{(j)}) \bigr\rangle \Bigr)  \, p_j(dx_j\mid x^{(j-1)}), \qquad (s \in {\Bbb R}^{m}, \, 1 \leq k < j \leq n)
$$
where as above $\partial / \partial x_k $ denotes the gradient with respect to $x_k \in {\Bbb R}^{m}$. 
The main result in this section is the following.

\begin{theorem}\label{thm5.1}
Suppose that there exist constants $\alpha >0$ and
$\kappa_{j,k} \geq0$ for $1\leq k<j\leq n$ such that

(i) $P^{(1)}$ and $p_k(. \mid x^{(k-1)})$ ($k=2, \dots, n; \, x^{(k-1)} \in {\Bbb R}^{m(k-1)}$) satisfy 
$LSI(\alpha )$;

(ii) $\Lambda_{j,k}(s) \leq \exp ( \kappa_{j,k} \, \Vert s \Vert^2/2)$ holds for all $s\in {\Bbb R}^m$.

\noindent Then the joint distribution $P^{(n)}$ satisfies
$LSI(\alpha_n)$ with 
\begin{equation}\label{5.1}
\alpha_n={{\alpha}\over {1+\varepsilon}} \Bigl(1+\sum_{k=0}^{n-2}
\prod_{m=k+1}^{n-1}(1+K_m)\Bigr)^{-1}
\end{equation}
\for all $\varepsilon >0,$ where $K_j={(1+\varepsilon^{-1})}\sum_{\ell =0}^{j-1}\kappa_{n-\ell ,
n-j}/\alpha $ for $j=1, \dots , n-1$.

Suppose further that there exist $R\geq 0$ and $\rho_\ell\geq 0$ for $\ell =1,
\dots , n-1$ such that 

\indent (iii) $\kappa_{j,k}\leq \rho_{j-k}$ for $1\leq k<j\leq n$, and  
$\sum_{\ell =1}^{n-1}\sqrt{\rho_\ell} \leq\sqrt R$.

\noindent Then $P^{(n)}$ satisfies $LSI(\alpha_n)$ where
$$
\alpha_n = \left\{ \begin{array}{ll}
		 \bigl( \sqrt{\alpha} -\sqrt {R}\bigr)^2  & \qquad \mbox{if} \quad R < \alpha,\\
	\displaystyle  \frac{\alpha}{n(n+1)(e-1)}         & \qquad \mbox{if} \quad R = \alpha,\\  
	\displaystyle	\Bigl(\frac{\alpha}{R}\Bigr)^n \frac{R -\alpha}{e(n+1)}
						    	  & \qquad \mbox{if}  \quad R > \alpha.
			\end{array}
			\right.			    
$$
\end{theorem}

\smallskip

Before proving this theorem, we give simple sufficient conditions for hypothesis 
(ii) to hold. When $m=1$, hypothesis (i) is equivalent to a condition on the cumulative
distribution functions by the criterion for $LSI$ given in \cite{BG99}.

\begin{proposition}\label{prop5.2}
In the above notation, let $1 \leq k <j$ and suppose that there exist
$\alpha>0$ and $L_{j,k}\geq 0$ such that

(i) $p_j(. \mid x^{(j-1)})$ satisfies $GC(1/\alpha )$ for all $x^{(j-1)} \in {\Bbb R}^{m(j-1)}$;

(ii) $u_j$ is twice continuously differentiable and the off-diagonal blocks of
its Hessian matrix satisfy  
$$ \Bigl\Vert{{\partial^2 u_{j}}\over {\partial x_j\partial x_k}}\Bigr\Vert \leq L_{j,k}
$$ 
as matrices $({\Bbb R}^m, \ell^2)\rightarrow ({\Bbb R}^m, \ell^2)$.

\noindent Then
$$
\Lambda_{j,k}(s)\leq \exp({{L^2_{j,k}} \, \Vert s \Vert^2/({2\alpha})})\qquad (s\in {\Bbb R}^m).
$$ 
\end{proposition}

\begin{proof}[Proof of Proposition \ref{prop5.2}]
Letting $s = \Vert s \Vert \, e$ for 
some unit vector $e$, we note that by (ii) 
the real function $x_j\mapsto \langle e , {{\partial u_j}/{\partial x_k}} \rangle $ is $L_{j,k}$-Lipschitz in the
variable of integration, and that  
$$
\int_{\Bbb R}\bigl\langle e , {{\partial u_j}\over {\partial x_k}}\bigr\rangle \, p_j(dx_j\mid x^{(j-1)})
= -\bigl\langle e , {{\partial}\over {\partial x_k}} \int_{\Bbb R}p_j(dx_j\mid x^{(j-1)}) \bigr\rangle=0
$$
\noindent since $p_j(. \mid x^{(j-1)})$ is a probability measure.
Then, by (i), \par
$$
\int_{\Bbb R}\exp \Bigl( \bigl\langle s , {{\partial u_j}\over {\partial x_k}} \bigr\rangle \Bigr) \, p_j(dx_j\mid x^{(j-1)}) 
\leq \exp ( \kappa \, L_{j,k}^2 \, \Vert s \Vert^2/2)
$$ 
\noindent holds for all $x^{(j-1)}$ in ${\Bbb R}^{m(j-1)}$. This inequality implies the Proposition.
\end{proof}

\smallskip

\begin{proof}[Proof of Proposition \ref{prop5.2}]
For notational convenience, $X$ denotes the 
state space ${\Bbb R}^{m}$. Then let 
$f:X^n\rightarrow {\Bbb R}$ be a smooth and compactly
supported function, and let $g_j:X^{n-j}\rightarrow {\Bbb R}$ be defined
by $g_0=f$ and by  
\begin{equation}\label{5.2}
g_j(x^{(n-j)})=\Bigl( \int_Xg_{j-1}(x^{(n-j+1)})^2 \, p_{n-j}(dx_{n-j+1}\mid
x^{(n-j)})\Bigr)^{1/2}
\end{equation}
for $j = 1, \dots, n-1$; finally, let $g_n$ be the constant $(\int f^2 dP^{(n)})^{1/2}$. \par
\indent From the recursive formula \eqref{5.2} one can easily verify the
identity 
\begin{equation}\label{5.3}
\int_{X^{n}}f^2\log \Bigl( f^2/\int_{X^{n}}f^2 dP^{(n)}\Bigr) dP^{(n)}
=\sum_{j=0}^{n-1}\int_{X^{n-j}}g_j^2\log 
\bigl( g_j^2/g_{j+1}^2\bigr) dP^{(n-j)}
\end{equation}
which is crucial to the proof; indeed, it allows us to obtain the result from
logarithmic Sobolev inequalities on $X$.

By hypothesis (i), the measure $dp_{n-j}=p_{n-j}(dx_{n-j}\mid x^{(n-j-1)})$ satisfies
$LSI(\alpha )$, whence
\begin{equation}\label{5.4}
\int_X g_j^2 \log \bigl( g_j^2/g_{j+1}^2\bigr) dp_{n-j}\leq 
{{2}\over {\alpha}}\int_X
\Bigl({{\partial g_j}\over {\partial x_{n-j}}}\Bigr)^2 dp_{n-j}\qquad 
(j=0, \dots, n-1),
\end{equation}
where for $j=n-1$ we take $dp_1=P^{(1)}(dx_1)$. The next step is to express these derivatives 
in terms of the gradient of $f$, using the identity
\begin{equation}\label{5.5}
g_j{{\partial g_j}\over {\partial x_{n-j}}}=\int_{X^{n-j}} f{{\partial f}\over
{\partial x_{n-j}}} \, dp_n\dots dp_{n-j+1} -{{1}\over {2}}\sum_{\ell
=0}^{j-1}\int_{X^{j-\ell}} g_\ell^2 \, {{\partial u_{n-\ell}}\over {\partial x_{n-j}}} \, dp_{n-\ell}\dots
dp_{n-j+1}
\end{equation}
which follows from the definition \eqref{5.2} of $g^2_j$ and that of $p_{n-j}$.
The integrals on the right-hand side of \eqref{5.5} will be bounded by the following Lemma.

\smallskip

\begin{lemma}\label{lem5.3}
Let $0 \leq \ell < j \leq n-1$, and assume that hypothesis (ii) holds. Then
\begin{equation}\label{5.6}
\Bigl\Vert\int_{X} g_\ell^2 \, {{\partial u_{n-\ell}}\over {\partial x_{n-j}}} \, dp_{n-\ell}\Bigr\Vert
\leq g_{\ell+1}\Bigl( 2\, \kappa_{n-\ell,n-j} \int_{X} g_\ell^2 \log (g_\ell^2/g_{\ell +1}^2)
\, dp_{n-\ell}\Bigr)^{1/2}.
\end{equation}
\end{lemma}

\begin{proof}[Proof of Lemma \ref{lem5.3}] 
By definition of $\Lambda_{n-\ell,n-j}$, we have
$$
\int_{X}\exp\Bigl( \bigl\langle s , {{\partial u_{n-\ell}}\over {\partial
x_{n-j}}} \bigr\rangle -\log \Lambda_{n-\ell,n-j}(s)\Bigr) dp_{n-\ell}\leq 1\qquad 
(s \in X),
$$
and hence by the dual formula for relative entropy, as in \cite[p. 693]{BGL01},
$$
\int_{X}\Bigl( \bigl\langle s , {{\partial u_{n-\ell}}\over {\partial x_{n-j}}}  \bigr\rangle 
-\log \Lambda_{n-\ell ,n-j}(s) \Bigr) g_\ell^2 \, dp_{n-\ell}
\leq \int_{X}{g_\ell^2}\log \bigl({{g_\ell^2}/{g_{\ell+1}^2}}\bigr)dp_{n-\ell}.
$$
Then hypothesis (ii) of the Theorem ensures that
$$
\bigl\langle s , \int_{X} {{\partial u_{n-\ell}}\over{\partial x_{n-j}}}
\, g_\ell^2 \, dp_{n-\ell} \bigr\rangle
\leq {{\Vert s \Vert^2}\over {2}} \, {\kappa_{n-\ell, n-j} \, g^2_{\ell+1}} + \int_X 
{g_\ell^2}\log\bigl({{g_\ell^2}/{g_{\ell+1}^2}}\bigr) dp_{n-\ell}
$$
and the stated result follows by optimizing this over $s \in {\Bbb R}^m$.
\end{proof}

\smallskip
\indent {\sl Conclusion of the Proof of Theorem 5.1.} When we integrate \eqref{5.6} with 
respect to $dp_{n-\ell -1}\dots dp_{n-j+1},$ we deduce
by the Cauchy--Schwarz inequality that
\begin{multline*}
\Bigl\vert\int_{{\Bbb R}^{j-\ell}} g_\ell^2{{\partial
u_{n-\ell}}\over{\partial x_{n-j}}}dp_{n-\ell}\dots dp_{n-j+1}\Bigr\vert \\
 \leq g_{j}\Bigl( 2 \, \kappa_{n-\ell,n-j}
\int_{{\Bbb R}^{j-\ell}} g_\ell^2 \log (g_\ell^2/g_{\ell +1}^2) dp_{n-\ell}\dots 
dp_{n-j+1}\Bigr)^{1/2}.
\end{multline*}
Then, by integrating the square of \eqref{5.5} with respect to $dP^{(n-j)}$ and making a further 
application of the Cauchy--Schwarz inequality, we obtain
\begin{equation}\label{5.7}
\int_{X^{n-j}}\Bigl\Vert {{\partial g_j}\over {\partial x_{n-j}}} \Bigr\Vert^2 dP^{(n-j)} 
\leq (1+\varepsilon) \! \int_{X^{n}}\Bigl\Vert {{\partial f}\over {\partial x_{n-j}}}\Bigr\Vert^2 dP^{(n)}
+{{1+\varepsilon^{-1}}\over {4}}\biggl\{ \sum_{\ell =0}^{j-1} 
\Bigl(2 \, \kappa_{n-\ell, n-j} \, h_{\ell}\Bigr)^{1/2}\biggr\}^{2}  
\end{equation}
where $\varepsilon >0$ is arbitrary and $h_{\ell}$ is given by
$$
h_{\ell} = \int_{X^{n-\ell}} g_\ell^2 \log \bigl({{g_\ell^2}/{g^2_{\ell +1}}}\bigr) dP^{(n-\ell )}.
$$

From \eqref{5.7}, which holds true for $j=1, \dots, n-1$, we first
prove the general result given in \eqref{5.1}. By \eqref{5.4} and the Cauchy--Schwarz inequality again, we obtain 
from \eqref{5.7} the crucial inequality
$$
h_j\leq d_j+K_j\sum_{m=0}^{j-1} h_m\qquad (j=1, \dots, n-1)
$$
where we have let
\begin{eqnarray*}
d_j&= &{{2(1+\varepsilon )}\over {\alpha }}\int_{{\Bbb R}^{n}}\Bigl( {{\partial f}\over 
{\partial x_{n-j}}}\Bigr)^2 dP^{(n)}\qquad (j=0,\dots, n-1), \\
K_j &=&{{1+\varepsilon^{-1}}\over {\alpha}}\sum_{\ell =0}^{j-1}\kappa_{n-\ell ,
n-j}\qquad\qquad (j=1,\dots, n-1).
\end{eqnarray*}
Since $h_0\leq d_0$ and all terms are positive, the partial
sums $H_k=\sum_{j =0}^k h_j $ satisfy the system of inequalities 
$$
H_k\leq d_k+(1+K_k)H_{k-1}\qquad (k=1,\dots, n-1),$$
with $H_0\leq d_0$. By induction, one can deduce that 
$$
H_{n-1}\leq d_{n-1}+\sum_{k=0}^{n-2}d_k\prod_{\ell =k+1}^{n-1}(1+K_\ell),
$$
 which in turn implies the bound
$$
H_{n-1}\leq \Bigl( 1+\sum_{k=0}^{n-2}\prod_{\ell =k+1}^{n-1}(1+K_\ell )\Bigr)
\sum_{j=0}^{n-1}d_j.
$$
By \eqref{5.3} this is equivalent to the inequality
$$
\int_{X^n}f^2\log \Bigl( f^2/\int_{X^n}f^2dP^{(n)}\Bigr) dP^{(n)}\leq 
{{2(1+\varepsilon )}\over {\alpha }}\Bigl(1+\sum_{k=0}^{n-2}\prod_{\ell =k+1}^{n-1}(1+K_\ell )\Bigr)
\int_{X^n}\Vert \nabla f\Vert^2 dP^{(n)}.
$$
Since $f$ is arbitrary, this ensures that $P^{(n)}$ satisfies
$LSI(\alpha_n)$ with $\alpha_n$ as in \eqref{5.1}.\par

\smallskip

\indent (iii) The extra hypothesis (iii) enables us to strengthen the preceding 
inequalities, so \eqref{5.7} leads to the convolution-type inequality
$$
h_j\leq d_j+{{1+\varepsilon^{-1}}\over {\alpha}}\Bigl( \sum_{\ell
=0}^{j-1}\sqrt{\rho_{j-\ell}}\,\sqrt {h_\ell}\Bigr)^2
$$
for $j=1, \dots, n-1$, and $h_0 \leq d_0$ for $j=0$. By summing over $j$ we obtain
$$ 
\sum_{j=0}^{k}h_j\leq \sum_{j=0}^{k}d_j+
{{1+\varepsilon^{-1}}\over {\alpha}}\sum_{j=1}^{k}\Bigl( \sum_{\ell
=0}^{j-1}\sqrt{\rho_{j-\ell}}\,\sqrt {h_\ell}\Bigr)^2,
$$
which implies by Young's convolution inequality that 
$$ 
\sum_{j=0}^{k}h_j\leq \sum_{j=0}^{k}d_j+
{{1+\varepsilon^{-1}}\over {\alpha}}\Bigl( \sum_{\ell
=1}^{k}\sqrt{\rho_{\ell}}\Bigr)^2 \sum_{\ell =0}^{k-1}{h_\ell}.
$$
Now let $R_j=(\sum_{\ell =1}^{j}\sqrt{\rho_\ell})^2$ 
 and $D_j=\sum_{\ell =0}^j d_\ell$; then by induction one can prove that 
$$
H_{k}\leq D_k+\sum_{j=0}^{k-1}D_j\prod_{\ell =j+1}^k {{1+\varepsilon^{-1}}\over
{\alpha}} R_{\ell}
$$
for $k=1, \dots, n-1$, and hence 
\begin{equation}\label{5.8}
H_{n-1}\leq \biggl( 1+\sum_{j=0}^{n-2}\Bigl(
{{1+\varepsilon^{-1}}\over {\alpha}} R\Bigr)^{n-j-1}\biggr) D_{n-1} =
\sum_{\ell=0}^{n-1} \Bigl(
{{1+\varepsilon^{-1}}\over {\alpha}} R\Bigr)^{\ell}  D_{n-1}
\end{equation}
since $D_j\leq D_{n-1}$ and $R_j\leq R$ by hypothesis
(iii). We finally select $\varepsilon$ to make the bound \eqref{5.8} precise, 
according to the relative values of $R$ and $\alpha$.

When $R=0$, we recover $LSI(\alpha )$ for $P^{(n)}$ as expected, since 
here $P^{(n)}$ is the tensor product of its marginal distributions, which satisfy
$LSI(\alpha )$.

When $0<R<\alpha$, we choose $\varepsilon =(\sqrt{(\alpha /R)}-1)^{-1}>0$ so that
$(1+\varepsilon^{-1})R/\alpha =\sqrt {({R}/{\alpha})}<1$ and hence
$$
H_{n-1}\leq D_{n-1}\sum_{\ell =0}^\infty (R/\alpha )^{\ell /2}={{D_{n-1}}\over
{1-\sqrt {(R/\alpha)}}},
$$
\noindent which by \eqref{5.3} and the definition of $H_{n-1}$ and $D_{n-1}$ implies the inequality 
$$
\int_{X^n}  f^2\log \Bigl( f^2/\int_{X^n}f^2dP^{(n)}\Bigr) dP^{(n)}\leq {{2}\over
{(\sqrt{\alpha}-\sqrt{R})^2}}\int_{X^n} \Vert \nabla f\Vert^2
dP^{(n)}.
$$

When $R\geq \alpha $, we choose $\varepsilon =n$ in \eqref{5.8}, obtaining
$$
H_{n-1}\leq {{2(n+1)}\over {\alpha}}
\biggl( {{(1+1/n)^n(R/\alpha )^n-1}\over {(1+1/n)(R/\alpha)-1}}\biggr)
\int_{X^n} \Vert \nabla f\Vert^2 dP^{(n)};
$$
as above this leads to the stated result by \eqref{5.3}. This concludes the proof.
\end{proof}

\medskip

\section{Logarithmic Sobolev inequalities for Markov processes}\label{sec6}
The results of the preceding section simplify considerably when we have an
homogeneous Markov process $(\xi_j)_{j=1}^n$ with state space ${\Bbb R}^m$, as we shall
now show. Suppose that the transition measure is $p(dy\mid x)=e^{-u(x,y)}dy$ where $u$ is a twice 
continuously differentiable function such that 
\begin{equation}\label{6.1}
\Lambda (s \mid x)=\int_{\Bbb R} \exp\Bigl( \bigl\langle s , {{\partial u}\over {\partial x}}
(x,y) \bigr\rangle \Bigr) p(dy\mid x)<\infty\qquad (s, x\in {\Bbb R}^m).
\end{equation}
Then Theorem  \ref{thm5.1} has the following consequence.

\smallskip

\begin{corollary}\label{cor6.1} 
Suppose that there exist constants $\kappa\geq 0$  and
$\alpha >0$ such
that:

(i) $P^{(1)}$ and $p(. \mid x)$ $(x\in {\Bbb R}^m)$ satisfy $LSI(\alpha )$;

(ii) $\Lambda (s\mid x)\leq \exp ({\kappa \Vert s \Vert^2/2})$ holds for
 all $s,x\in {\Bbb R}^m$.
 
\noindent Then the joint law $P^{(n)}$ of the first $n$ variables satisfies $LSI(\alpha_n)$,
 where 

$$
\alpha_n = \left\{ \begin{array}{ll}
		 \bigl( \sqrt{\alpha} -\sqrt {\kappa}\bigr)^2  & \qquad \mbox{if} \quad \kappa < \alpha,\\
	\displaystyle  \frac{\alpha}{n(n+1)(e-1)} 	       & \qquad \mbox{if} \quad \kappa = \alpha,\\  
	\displaystyle	\Bigl(\frac{\alpha}{\kappa}\Bigr)^n \frac{\kappa -\alpha}{e(n+1)}
						    	       & \qquad \mbox{if}  \quad \kappa > \alpha.
			\end{array}
			\right.			    
$$
\end{corollary}

\smallskip

\begin{proof}
In the notation of section 5, we have $u_j(x^{(j)})=u(x_{j-1}, x_j)$, so we can
take $\kappa_{j,m}=0$ for $m=1, \dots, j-2$, and $\kappa_{j,j-1}=\kappa$ for
 $j=2, \dots, n$; hence we can take $\rho_1=\kappa$ and $\rho_j=0$ for $j=2, 3, \dots.$  
 Now we can apply Theorem \ref{thm5.1} (iii) and obtain the stated result with $R=\kappa$ 
in the various cases. (In fact \eqref{5.7} simplifies considerably for a Markov process,
and hence one can obtain an easier direct proof of Corollary \ref{cor6.1}.)
\end{proof}

\smallskip

\begin{proof}[Proof of Theorem \ref{thm1.3}]

By the mean-value theorem and hypothesis (ii) of Theorem \ref{thm1.3}, the function $y\mapsto \langle e , {\partial u}/{\partial x} \rangle$ is 
$L$-Lipschitz $({\Bbb R}^m, \ell^2)\rightarrow {\Bbb R}$ for any unit vector
$e$ in ${\Bbb R}^m$, and hence $\Lambda (s\mid x)\leq \exp
( \Vert s \Vert^2L^2/(2\alpha ))$ holds for all $s\in {\Bbb R}^m$ as in Proposition \ref{prop5.2}. Hence we can
take $\kappa =L^2/\alpha $ in Corollary \ref{cor6.1} and deduce Theorem \ref{thm1.3} with the various
values of the constant.
\end{proof}

\smallskip

\begin{remarks}\label{rem6.2} 
(i) Theorem \ref{thm5.1} and Corollary \ref{cor6.1} extend with
suitable changes in notation when the state space is a connected $C^1$-smooth 
Riemannian manifold $X$. The proofs reduce to calculations in local co-ordinate
charts. McCann \cite{MC01} has shown that a locally Lipschitz function on $X$ is differentiable
except on a set that has zero Riemannian volume; so a $L$-Lipschitz condition on 
$f:X\rightarrow {\Bbb R}$ is essentially equivalent
to $\Vert \nabla f\Vert\leq L.$

(ii) Corollary 6.1 is a natural refinement of Theorems \ref{thm1.1} and \ref{thm1.2}. Indeed $LSI(\alpha)$ implies 
$T_s(\alpha)$. Then, in the notation of the mentioned results, suppose that $u$ is a twice 
continuously differentiable function with bounded second-order partial derivatives. Then, by Proposition
\ref{prop2.2}, hypotheses (i) and (ii) of Corollary \ref{cor6.1} together imply that the map $x\mapsto p(.\mid x)$ 
is $(\kappa/\alpha )^{1/2}$ Lipschitz as a function ${\Bbb R}^m \rightarrow ({\hbox{Prob}}_2\,({\Bbb R}^m), W_2)$,
hence ${\Bbb R}^m\rightarrow ({\hbox{Prob}}_s\,({\Bbb R}^m), W_s)$ as in Theorems \ref{thm1.1} or \ref{thm1.2}. 
Similarly Proposition \ref{prop2.2} ensures that
Theorem \ref{thm5.1} is a refinement of Theorem \ref{thm2.1} with, for $s=2$, 
$$
M \leq M_2/\alpha ={{1}\over {\alpha}}
\sup_{x^{(j-1)}}\sum_{k=1}^{j-1}\int_{{\Bbb R}^m}\Bigl\Vert{{\partial u_j}\over{\partial x_k}}\Bigr\Vert^2 \,
p_j(dx_j\mid x^{(j-1)})\leq{{1}\over
{\alpha}}\sum_{k=1}^{j-1}\kappa_{j,k}.
$$
Note also the similarity between the constants in Theorem \ref{thm2.1} (iii) and
Theorem \ref{thm5.1} (iii) when $s=2$ and one rescales $R$ suitably. In Example \ref{ex6.3} we show these constants 
to be optimal.
\end{remarks}

\smallskip

\begin{example}\label{ex6.3} 
(Ornstein--Uhlenbeck Process) 
We now show that the constants $\kappa_n$ of Theorem \ref{thm1.1} (or Theorem \ref{thm3.1}(iii)) and $\alpha_n$ of
Corollary \ref{cor6.1} have optimal growth in $n$. For this purpose we consider the real 
Ornstein--Uhlenbeck process conditioned to
start at $x \in {\Bbb R}$, namely the solution to the It\^o stochastic differential equation
$$ dZ_t^{(x)} =-\rho Z_t^{(x)} dt+dB_t^{(0)}, \qquad (t \geq 0)$$
where $(B_t^{(0)})$ is a real standard Brownian motion starting at $0$, and 
$\rho\in {\Bbb R}$. In financial modelling, OU processes with $\rho <0$
are used to model stock prices in a rising market (see \cite[p 26]{Fol95} for instance). More precisely 
we consider the discrete-time Markov process $(\xi_j)_{j=1}^n$ defined by 
$\xi_j= Z_{j\tau}^{(x)}$ where $\tau>0$, and test the Gaussian
concentration inequality with 
the $1$-Lipschitz function $F_n:(\Bbb{R}^n, \ell^1) \rightarrow {\Bbb R}$ defined by 
$F_n(x^{(n)}) = \sum_{j=1}^{n} x_j$. \par
\indent The exponential integral satisfies
\begin{equation}\label{6.2}
\int_{{\Bbb R}^n} \exp\bigl({s F_n(x^{(n)})}\bigr) \, P^{(n)}(dx^{(n)}) = 
\Bbb{E} \, \exp\bigl({s F_n(\xi^{(n)})}\bigr) =
\Bbb{E} \,\exp\Bigl({s \sum_{j=1}^{n} Z_{j\tau}^{(x)}}\Bigr).
\end{equation}
This sum can be expressed in terms of the increments of the OU process
$$
\sum_{j=1}^{n} Z_{j\tau}^{(x)} = \sum_{i=1}^{n} \theta^i \, Z^{(x)}_0  + \sum_{j=0}^{n-1} 
\, \sum_{i=0}^{n-j-1} \theta^i  \, \bigl(Z_{(j+1)\tau}^{(x)} - \theta Z_{j\tau}^{(x)}\bigr),
$$
\noindent with $\theta = e^{-\rho \tau}$. Moreover one can integrate the stochastic
differential equation and prove that 
$(Z_{(j+1)\tau}^{(x)} - \theta Z_{j\tau}^{(x)})_{0\leq k \leq n-1}$ 
are independent random variables each with $N(0,\sigma^2)$ distribution, where 
$\sigma^2 =({1-\theta^2})/({2 \, \rho})$ when $\rho \neq 0$, and $\sigma^2=\tau$ when $\rho =0$.
Hence the exponential integral \eqref{6.2} equals 
$$ 
\exp\Bigl(s \sum_{i=1}^{n}\theta^i x\Bigr)\prod_{j=0}^{n-1} 
{\Bbb E}\exp \Bigl[s\Bigl(\sum_{i=0}^{n-j-1} \theta^i \Bigr)
\bigl(Z_{(j+1)\tau}^{(x)} - \theta Z_{j\tau}^{(x)}\bigr)  \Bigr] 
= \exp \bigl(s \, \Bbb{E} \, F_n(\xi^{(n)}) + s^2 \kappa_n/2 \bigr)
$$
where
\begin{equation}\label{6.3}
\kappa_n = \sigma^2\sum_{j=0}^{n-1} \Bigl(\sum_{i=0}^{n-j-1} \theta^i \Bigr)^2.
\end{equation}

However, hypothesis (i) of Theorem \ref{thm1.1} holds with $L=\theta$, since 
$P^{(1)}$ with distribution $N(x,\sigma^2)$ and $p(. \mid x)$ with 
distribution $N(\theta x,\sigma^2)$ satisfy
$GC(\kappa_1)$ where $\kappa_1 = \sigma^2$, while hypothesis (ii) is satisfied with  
\begin{equation}\label{6.4}
W_1(p(. \mid x) , p(. \mid x')) = 
W_1(N(\theta x,\sigma^2), N(\theta x',\sigma^2)) = \theta \vert x-x'
\vert \qquad (x,x'\in {\Bbb R}).
\end{equation}
Hence the constant $\kappa_n(L)$ given by Theorem 1.1 is exactly
the directly computed constant $\kappa_n$ in \eqref{6.3}, in each of the cases $L=1$, $L>1$ and $L<1$, 
corresponding to $\rho =0$, $\rho <0$ and $\rho >0$.

\smallskip

As regards Corollary \ref{cor6.1}, note that the transition probability is given by
$$
p(dy\mid x)={{1}\over{\sqrt{2 \pi \sigma^2}}} \exp\Bigl( -{{(y-\theta x)^2}\over{2 \sigma^2}}\Bigr) dy
$$
since $Z_\tau^{(x)}$ is distributed as $\theta x+B^{(0)}_{\sigma^2}$. Hence by direct calculation
we have
$$
\alpha ={{1}\over {\sigma^2}}, \quad\kappa =
{{\theta^2}\over {\sigma^2}},\quad L={{\theta}\over{\sigma^2}};
$$
consequently the dependence parameters $(\kappa /\alpha )^{1/2}$ and $\theta$ given in \eqref{6.4}
coincide, as in Remark \ref{rem6.2}(ii). 

Further, by considering the function 
$f(x^{(n)}) =\exp \bigl(\sum_{j=1}^{n} \theta^j x_j \bigr)$, 
one can prove that the joint law $P^{(n)}$ cannot satisfy a logarithmic Sobolev 
inequality with $\alpha_n$ greater than
some constant multiple of $n^{-3}$ for $\theta = 1$, and $(\alpha / \kappa)^n$ 
for $\theta >1$. Thus for $\theta \geq 1$,
we recover the order of growth in $n$ of the constants given in Corollary \ref{cor6.1}; whereas 
for $\theta < 1$, the constant given in Corollary \ref{cor6.1} is independent of $n$.

The OU process does not satisfy the Doeblin condition $D_0$, as Rosenblatt observes; see 
\cite[p. 214]{Ros71}.
\end{example}

\bigskip

\indent {\bf Acknowledgements.} This research was supported in part by the European Network 
PHD, MCRN -511953. The authors thank Professors P.J. Diggle, M. Ledoux, K.
Marton and C. Villani for helpful conversations.

\end{document}